\documentclass{article}
\usepackage[utf8]{inputenc}
\usepackage{amsmath}
\usepackage{amsfonts}
\usepackage{float}
\usepackage{bbm}
\usepackage{makecell}
\usepackage{array}
\usepackage{changepage}
\newcolumntype{?}{!{\vrule width 1.5pt}}
\begin{document}
\title{Identities of the Fractional Fourier Transform \\ and the Versor Transform}
\author{Maurice Pierre, Class of 2021}
\date{April 4, 2021}
\maketitle
\hspace{103pt} Supervisor: Peter Woit
\begin{abstract}
We provide an introduction to the Fractional Fourier Transform $\mathcal{F}_{\theta}$ and draw a connection between it and the unit complex number $e^{i\theta}$. Motivated by this, we define an entirely new object associated with any unit quaternion $e^{i\xi_{1}}\cos\eta+e^{i\xi_{2}}j\sin\eta$, which we call the Versor Transform $\mathcal{V}_{(\xi_{1},\eta,\xi_{2})}$. This transform, which has both the Fourier and Laplace Transforms as special cases, encourages an alternate view of the relationship between them. We also derive several identities for both $\mathcal{F}_{\theta}$ and $\mathcal{V}_{(\xi_{1},\eta,\xi_{2})}$.

\end{abstract}
\null
\textit{Note: For all of the transforms described in this paper, we will assume the input function $f(x)$ is sufficiently nice to ensure that they are well-defined.}
\section{Linear Canonical Transform}
\subsection{Definition}
The Linear Canonical Transform (LCT) is defined for a matrix \\ $\mathrm{M}=\big[\begin{smallmatrix}
a & b \\
c & d
\end{smallmatrix}\big]\in\mathrm{SL}_{2}(\mathbb{C})$ ($a,b,c,d\in\mathbb{C}$ and $|\mathrm{M}|=1$) as follows:
$$[\mathcal{LCT}_{\mathrm{M}}\{f(x)\}](u)=
\begin{cases}
\displaystyle\sqrt{\frac{1}{ib}}\int_{-\infty}^{\infty}e^{\pi i\frac{1}{b}(ax^2+du^2)}f(x)e^{-2\pi i\frac{1}{b}xu}\ dx, & b\neq0 \\ \\
\sqrt{d}e^{\pi icdu^2}f(du), & b=0
\end{cases}\quad\cite{1}$$
This object generalizes several classical transforms, including the Fourier and Laplace Transforms, as we will see later.
\newpage
\subsection{Identities}
The proofs of these identities are given in the appendix. Assume $b\neq0$.
\subsubsection{Shift Identities}
For a shift in the original function's domain, $f(x)\rightarrow f(x-x_0)$:
$$[\mathcal{LCT}_{\mathrm{M}}\{f(x-x_{0})\}](u)=e^{-\pi iacx_{0}^2}e^{2\pi icx_{0}u}[\mathcal{LCT}_{\mathrm{M}}\{f(x)\}](u-ax_{0})$$
For a shift in the transformed function's domain, \\ $[\mathcal{LCT}_{\mathrm{M}}\{f(x)\}](u)\rightarrow[\mathcal{LCT}_{\mathrm{M}}\{f(x)\}](u-u_{0})$:
$$e^{\pi i\frac{1}{b}du_{0}^2}e^{-2\pi i\frac{1}{b}duu_{0}}[\mathcal{LCT}_{\mathrm{M}}\{e^{2\pi i\frac{1}{b}xu_{0}}f(x)\}](u)=[\mathcal{LCT}_{\mathrm{M}}\{f(x)\}](u-u_{0})$$
\subsubsection{Derivative Identities}
For the $n^{\text{th}}$ derivative of the original function, $f(x)\rightarrow\displaystyle\frac{d^n}{dx^n}f(x)$:
\begin{adjustwidth}{-44pt}{-44pt}
$$\bigg[\mathcal{LCT}_{\mathrm{M}}\bigg\{\frac{d^n}{dx^n}f(x)\bigg\}\bigg](u)=(-1)^n\sum_{k=0}^{n}(-2\pi ib^{-1}u)^k\binom{n}{k}\bigg[\mathcal{LCT}_{\mathrm{M}}\bigg\{e^{-\pi i\frac{1}{b}ax^2}\frac{d^{n-k}}{dx^{n-k}}(e^{\pi i\frac{1}{b}ax^2})f(x)\bigg\}\bigg](u)$$
\end{adjustwidth}
For the $n^{\text{th}}$ derivative of the transformed function, \\ $[\mathcal{LCT}_{\mathrm{M}}\{f(x)\}](u)\rightarrow\displaystyle\frac{d^n}{du^n}\Big([\mathcal{LCT}_{\mathrm{M}}\{f(x)\}](u)\Big)$:
\begin{adjustwidth}{-27pt}{-27pt}
$$\sum_{k=0}^{n}\binom{n}{k}e^{-\pi i\frac{1}{b}du^2}\frac{d^{n-k}}{du^{n-k}}(e^{\pi i\frac{1}{b}du^2})[\mathcal{LCT}_{\mathrm{M}}\{(-2\pi ib^{-1}x)^kf(x)\}](u)=\frac{d^n}{du^n}\Big([\mathcal{LCT}_{\mathrm{M}}\{f(x)\}](u)\Big)$$
\end{adjustwidth}
The function $e^{-\pi i\frac{1}{b}ax^2}\frac{d^m}{dx^m}(e^{\pi i\frac{1}{b}ax^2})$ can be expressed in terms of Hermite polynomials as follows:
\begin{adjustwidth}{-5pt}{-5pt}
$$e^{-\pi i\frac{1}{b}ax^2}\frac{d^m}{dx^m}(e^{\pi i\frac{1}{b}ax^2})=(-1)^{m}He_{m}\Big(\sqrt{-2\pi iab^{-1}}x\Big)=(-1)^{m}H_{m}\Big(\sqrt{-\pi iab^{-1}}x\Big)$$
\end{adjustwidth}
where $He_m$ are the probabilist's version and $H_m$ are the physicist's version.
\subsubsection{Parseval's Theorem}
For any two functions $f(x)$ and $g(x)$ for which their LCTs exist:
$$\int_{-\infty}^{\infty}f(x)\overline{g(x)}\ dx=\int_{-\infty}^{\infty}[\mathcal{LCT}_{\mathrm{M}}\{f(x)\}](u)\ \overline{[\mathcal{LCT}_{\mathrm{M}}\{g(x)\}](u)}\ du$$
Letting $f(x)=g(x)$ gives Plancherel's Theorem:
$$\int_{-\infty}^{\infty}|f(x)|^2\ dx=\int_{-\infty}^{\infty}\Big|[\mathcal{LCT}_{\mathrm{M}}\{f(x)\}](u)\Big|^2\ du$$
\subsubsection{Convolution Theorem}
Let the convolution of two functions $f(x)$ and $g(x)$ be defined as follows:
$$[f*g](x)=\int_{-\infty}^{\infty}f(\tau)g(x-\tau)\ d\tau\text{. Then}$$
\begin{adjustwidth}{-1pt}{-1pt}
$$\sqrt{\frac{1}{ib}}\ e^{\pi i\frac{1}{b}du^2}[\mathcal{LCT}_{\mathrm{M}}\{e^{-\pi i\frac{1}{b}ax^2}[fh*gh](x)\}](u)=[\mathcal{LCT}_{\mathrm{M}}\{f(x)\}\cdot\mathcal{LCT}_{\mathrm{M}}\{g(x)\}](u)$$
\end{adjustwidth}
\begin{adjustwidth}{-9pt}{-9pt}
$$\sqrt{\frac{1}{-ib}}^{-1}e^{-\pi i\frac{1}{b}du^2}[\mathcal{LCT}_{\mathrm{M}}\{e^{\pi i\frac{1}{b}ax^2}[f\cdot g](x)\}](u)=[\mathcal{LCT}_{\mathrm{M}}\{f(x)\}k*\mathcal{LCT}_{\mathrm{M}}\{g(x)\}k](u)$$
\end{adjustwidth}
where $h(x)=e^{\pi i\frac{1}{b}ax^2}$ and $k(u)=e^{-\pi i\frac{1}{b}du^2}$.
\subsubsection{Cross-Correlation Theorem}
Let the cross-correlation of two functions $f(x)$ and $g(x)$ be defined as follows:
$$[f\star g](x)=\int_{-\infty}^{\infty}\overline{f(\tau)}g(x+\tau)\ d\tau\text{. Then}$$
\begin{adjustwidth}{-4pt}{-4pt}
$$\overline{\sqrt{\frac{1}{ib}}}\ e^{-\pi i\frac{1}{b}du^2}[\mathcal{LCT}_{\mathrm{M}}\{e^{-\pi i\frac{1}{b}ax^2}[fh\star gh](x)\}](u)=[\overline{\mathcal{LCT}_{\mathrm{M}}\{f(x)\}}\cdot\mathcal{LCT}_{\mathrm{M}}\{g(x)\}](u)$$
\end{adjustwidth}
\begin{adjustwidth}{-12pt}{-12pt}
$$\overline{\sqrt{\frac{1}{-ib}}}^{-1}e^{-\pi i\frac{1}{b}du^2}[\mathcal{LCT}_{\mathrm{M}}\{e^{-\pi i\frac{1}{b}ax^2}[\overline{f}\cdot g](x)\}](u)=[\mathcal{LCT}_{\mathrm{M}}\{f(x)\}k\star\mathcal{LCT}_{\mathrm{M}}\{g(x)\}k](u)$$
\end{adjustwidth}
where $h(x)=e^{\pi i\frac{1}{b}ax^2}$ and $k(u)=e^{-\pi i\frac{1}{b}du^2}$.
\section{Fractional Fourier Transform}
\subsection{Motivation}
The Fractional Fourier Transform has existed in the mathematical literature since the 1930s.\;\cite{2} This section is designed to be an introduction to the topic in the context of the standard Fourier Transform and the Linear Canonical Transform. \\\null\\
The standard Fourier Transform is defined as follows:
$$[\mathcal{F}\{f(x)\}](u)=\int_{-\infty}^{\infty}f(x)e^{-2\pi ixu}\ dx$$
Furthermore, composing the Fourier Transform with itself yields
$$\mathcal{F}^2=\mathcal{P},\quad \mathcal{F}^3=\mathcal{F}^{-1},\quad \mathcal{F}^4=\mathit{1}\quad\cite{4}$$
where $\mathit{1}:f(x)\rightarrow f(u)$ is the identity operator, $\mathcal{P}:f(x)\rightarrow f(-u)$ is the parity operator, and $\mathcal{F}^{-1}$ is the inverse Fourier Transform:
$$[\mathcal{F}^{-1}\{f(x)\}](u)=\int_{-\infty}^{\infty}f(x)e^{2\pi ixu}\ dx$$
This bears a striking resemblance to the imaginary unit $i$:
\begin{table}[H]
\centering
\begin{tabular}{|c?c|c|c|c|c|}
\hline
$n$ & $0$ & $1$ & $2$ & $3$ & $4$ \\ \hline
$i^n$ & $1$ & $i$ & $-1$ & $-i$ & $1$ \\ \hline
$\mathcal{F}^n$ & $\mathit{1}$ & $\mathcal{F}$ & $\mathcal{P}$ & $\mathcal{F}^{-1}$ & $\mathit{1}$ \\ \hline
\end{tabular}
\end{table}
\noindent In other words, $\mathcal{F}^2=\mathcal{P}$ is the Fourier Transform analogue of $i^2=-1$. \\\null\\
The Fourier Transform can also be thought of as the LCT with the input matrix $\mathrm{F}=\big[\begin{smallmatrix}
0 & 1 \\
-1 & 0
\end{smallmatrix}\big]$ (up to a factor of $e^{-\frac{\pi}{4}i}$). This matrix corresponds to a rotation by 90\textdegree\ clockwise in $\mathbb{R}^2$. We can generalize this by using an arbitrary clockwise rotation matrix $\mathrm{F}_{\theta}=\big[\begin{smallmatrix}
\cos\theta & \sin\theta \\
-\sin\theta & \cos\theta
\end{smallmatrix}\big],\ \theta\in\mathbb{R}$. This gives us the following:
$$[\mathcal{LCT}_{\mathrm{F}_{\theta}}\{f(x)\}](u)=\sqrt{-i\csc(\theta)}\int_{-\infty}^{\infty}e^{\pi i\cot(\theta)(x^2+u^2)}f(x)e^{-2\pi i\csc(\theta)xu}\ dx$$
The Fractional Fourier Transform $\mathcal{F}_{\theta}$ is defined so that $\mathcal{F}_{n\pi/2}$ agrees with $\mathcal{F}^n$ for all $n\in\mathbb{Z}$. To do this, the factor of $\sqrt{-i\csc(\theta)}$ is changed to $\sqrt{1-i\cot(\theta)}$.
$$[\mathcal{F}_{\theta}\{f(x)\}](u)=\sqrt{1-i\cot(\theta)}\int_{-\infty}^{\infty}e^{\pi i\cot(\theta)(x^2+u^2)}f(x)e^{-2\pi i\csc(\theta)xu}\ dx\quad\cite{4}$$
The inverse transform $\mathcal{F}^{-1}_{\theta}$ is equivalent to $\mathcal{F}_{-\theta}$. This can be seen by the relationship between their underlying matrices:
$$\mathrm{F}^{-1}_{\theta}=\begin{bmatrix}
\cos\theta & -\sin\theta \\
\sin\theta & \cos\theta
\end{bmatrix}=
\begin{bmatrix}
\cos(-\theta) & \sin(-\theta) \\
-\sin(-\theta) & \cos(-\theta)
\end{bmatrix}=\mathrm{F}_{-\theta}$$
Furthermore, it is well-known that a matrix of the form $\big[\begin{smallmatrix}
a & b \\
-b & a
\end{smallmatrix}\big]$ with $a,b\in\mathbb{R}$ may be interpreted as the complex number $a+bi$. Doing this for $\mathrm{F}_{\theta}$ yields $\cos\theta+i\sin\theta=e^{i\theta}$. This means that every unit complex number $e^{i\theta}$ has an associated transform $\mathcal{F}_{\theta}$. \\\null\\
The standard Fourier Transform of a function takes it from the time domain into the frequency domain. As a result, the domain of the Fractional Fourier Transform can be considered a linear combination of both time and frequency. \\\null\\
It should be noted that when $\theta$ is an integer multiple of $\pi$, the trigonometric values in $\mathcal{F}_{\theta}$ diverge. However, this can be resolved by thinking of $\mathcal{F}_{\theta}$ as the integral of $f(x)$ multiplied by a kernel function $K_{\theta}(x,u)$:
$$[\mathcal{F}_{\theta}\{f(x)\}](u)=\int_{-\infty}^{\infty}K_{\theta}(x,u)f(x)\ dx$$
$$\text{where } K_{\theta}(x,u)=\sqrt{1-i\cot(\theta)}e^{\pi i\cot(\theta)(x^2+u^2)}e^{-2\pi i\csc(\theta)xu}.$$
By taking limits, one can see that
$$K_{\theta}(x,u)=
\begin{cases}
\delta(x-u), & \theta\equiv0\mod2\pi \\
\delta(x+u), & \theta\equiv\pi\mod2\pi
\end{cases}\quad\cite{3}$$
where $\delta$ is the Dirac delta function. Integrating $f(x)$ multiplied by these delta functions gives
$$[\mathcal{F}_{\theta}\{f(x)\}](u)=
\begin{cases}
[\mathit{1}\{f(x)\}](u), & \theta\equiv0\mod2\pi \\
[\mathcal{P}\{f(x)\}](u), & \theta\equiv\pi\mod2\pi
\end{cases}$$
where $\mathit{1}$ and $\mathcal{P}$ are the identity and parity operators mentioned previously.
\subsection{Identities}
These can be derived by letting
$$\mathrm{M}=\mathrm{F}_{\theta}=\begin{bmatrix}
\cos\theta & \sin\theta \\
-\sin\theta & \cos\theta
\end{bmatrix}$$
in the corresponding identities for the LCT. For the Convolution and Cross-Correlation Theorems, the constant factors are adjusted appropriately.
\subsubsection{Shift Identities}
For a shift in the original function's domain,
$$[\mathcal{F}_{\theta}\{f(x-x_{0})\}](u)=e^{\pi i\cos(\theta)\sin(\theta)x_{0}^2}e^{-2\pi i\sin(\theta)x_{0}u}[\mathcal{F}_{\theta}\{f(x)\}](u-\cos(\theta)x_{0})$$
For a shift in the transformed function's domain,
$$e^{\pi i\cot(\theta)u_{0}^2}e^{-2\pi i\cot(\theta)uu_{0}}[\mathcal{F}_{\theta}\{e^{2\pi i\csc(\theta)xu_{0}}f(x)\}](u)=[\mathcal{F}_{\theta}\{f(x)\}](u-u_{0})$$
\newpage
\subsubsection{Derivative Identities}
For the $n^{\text{th}}$ derivative of the original function,
\begin{adjustwidth}{-45pt}{-45pt}
$$\bigg[\mathcal{F}_{\theta}\bigg\{\frac{d^n}{dx^n}f(x)\bigg\}\bigg](u)=(-1)^n\sum_{k=0}^{n}(-2\pi i\csc(\theta)u)^k\binom{n}{k}\bigg[\mathcal{F}_{\theta}\bigg\{e^{-\pi i\cot(\theta)x^2}\frac{d^{n-k}}{dx^{n-k}}(e^{\pi i\cot(\theta)x^2})f(x)\bigg\}\bigg](u)$$
\end{adjustwidth}
For the $n^{\text{th}}$ derivative of the transformed function,
\begin{adjustwidth}{-28pt}{-28pt}
$$\sum_{k=0}^{n}\binom{n}{k}e^{-\pi i\cot(\theta)u^2}\frac{d^{n-k}}{du^{n-k}}(e^{\pi i\cot(\theta)u^2})[\mathcal{F}_{\theta}\{(-2\pi i\csc(\theta)x)^kf(x)\}](u)=\frac{d^n}{du^n}\Big([\mathcal{F}_{\theta}\{f(x)\}](u)\Big)$$
\end{adjustwidth}
\subsubsection{Parseval's Theorem}
For any two functions $f(x)$ and $g(x)$ for which their Fractional Fourier Transforms exist:
$$\int_{-\infty}^{\infty}f(x)\overline{g(x)}\ dx=\int_{-\infty}^{\infty}[\mathcal{F}_{\theta}\{f(x)\}](u)\ \overline{[\mathcal{F}_{\theta}\{g(x)\}](u)}\ du$$
Letting $f(x)=g(x)$ gives Plancherel's Theorem:
$$\int_{-\infty}^{\infty}|f(x)|^2\ dx=\int_{-\infty}^{\infty}\Big|[\mathcal{F}_{\theta}\{f(x)\}](u)\Big|^2\ du$$
\subsubsection{Convolution Theorem}
\begin{adjustwidth}{-6pt}{-6pt}
$$\sqrt{1-i\cot(\theta)}e^{\pi i\cot(\theta)u^2}[\mathcal{F}_{\theta}\{e^{-\pi i\cot(\theta)x^2}[fh*gh](x)\}](u)=[\mathcal{F}_{\theta}\{f(x)\}\cdot\mathcal{F}_{\theta}\{g(x)\}](u)$$
\end{adjustwidth}
\begin{adjustwidth}{-11pt}{-11pt}
$$\sqrt{1+i\cot(\theta)}^{-1}e^{-\pi i\cot(\theta)u^2}[\mathcal{F}_{\theta}\{e^{\pi i\cot(\theta)x^2}[f\cdot g](x)\}](u)=[\mathcal{F}_{\theta}\{f(x)\}k*\mathcal{F}_{\theta}\{g(x)\}k](u)$$
\end{adjustwidth}
where $h(x)=e^{\pi i\cot(\theta)x^2}$ and $k(u)=e^{-\pi i\cot(\theta)u^2}$.
\subsubsection{Cross-Correlation Theorem}
\begin{adjustwidth}{-10pt}{-10pt}
$$\overline{\sqrt{1-i\cot(\theta)}}\ e^{-\pi i\cot(\theta)u^2}[\mathcal{F}_{\theta}\{e^{-\pi i\cot(\theta)x^2}[fh\star gh](x)\}](u)=[\overline{\mathcal{F}_{\theta}\{f(x)\}}\cdot\mathcal{F}_{\theta}\{g(x)\}](u)$$
\end{adjustwidth}
\begin{adjustwidth}{-14pt}{-14pt}
$$\overline{\sqrt{1+i\cot(\theta)}}^{-1}e^{-\pi i\cot(\theta)u^2}[\mathcal{F}_{\theta}\{e^{-\pi i\cot(\theta)x^2}[\overline{f}\cdot g](x)\}](u)=[\mathcal{F}_{\theta}\{f(x)\}k\star\mathcal{F}_{\theta}\{g(x)\}k](u)$$
\end{adjustwidth}
where $h(x)=e^{\pi i\cot(\theta)x^2}$ and $k(u)=e^{-\pi i\cot(\theta)u^2}$.
\section{Versor Transform}
\subsection{Motivation}
The Versor Transform is an entirely new object defined for this paper. In this section, we motivate its definition by looking for a transform which corresponds to any unit quaternion. \\\null\\
It is known that the quaternion $a+bi+cj+dk$ with $a,b,c,d\in\mathbb{R}$ can be represented as a matrix of the form $\big[\begin{smallmatrix}
a+bi & c+di \\
-c+di & a-bi
\end{smallmatrix}\big]$ containing only complex entries. Setting one component to 1 and the rest to 0, we obtain the following four basis matrices:
$$\mathbbm{1}=\begin{bmatrix}
1 & 0 \\
0 & 1
\end{bmatrix},\quad
\mathrm{I}=\begin{bmatrix}
i & 0 \\
0 & -i
\end{bmatrix},\quad
\mathrm{J}=\begin{bmatrix}
0 & 1 \\
-1 & 0
\end{bmatrix},\quad
\mathrm{K}=\begin{bmatrix}
0 & i \\
i & 0
\end{bmatrix}$$
Since all of these matrices have determinant 1, each one has an associated LCT. For simplicity, we will ignore the constant factors of $\sqrt{1/ib}$ and $\sqrt{d}$.
$$[\mathit{1}\{f(x)\}](u)=f(u)$$
$$[\mathcal{I}\{f(x)\}](u)=f(u/i)$$
$$[\mathcal{J}\{f(x)\}](u)=\int_{-\infty}^{\infty}f(x)e^{-2\pi ixu}\ dx$$
$$[\mathcal{K}\{f(x)\}](u)=\int_{-\infty}^{\infty}f(x)e^{-2\pi xu}\ dx$$
These operators correspond to the four basis quaternions 1, $i$, $j$, and $k$. $\mathit{1}$ is the identity operator, $\mathcal{I}$ is a scaling of the domain by a factor of $i$, $\mathcal{J}$ is the Fourier Transform, and $\mathcal{K}$ is the Bilateral Laplace Transform scaled by a factor of $2\pi$. More specifically, $[\mathcal{K}\{f(x)\}](u)=[\mathcal{B}\{f(x)\}](2\pi u)$, where
$$[\mathcal{B}\{f(x)\}](u)=\int_{-\infty}^{\infty}f(x)e^{-xu}\ dx$$
Each individual operator $\mathcal{T}=\mathcal{I}$, $\mathcal{J}$, or $\mathcal{K}$ satisfies $\mathcal{T}^2=\mathcal{P}$, $\mathcal{T}^3=\mathcal{T}^{-1}$, and \\ $\mathcal{T}^4=\mathit{1}$, where $\mathcal{P}$ is the parity operator and $\mathcal{T}^{-1}$ is the inverse of $\mathcal{T}$. \\ Furthermore,
$$[\mathcal{I}\mathcal{J}\{f(x)\}](u)=\Big[\mathcal{I}\Big\{[\mathcal{J}\{f(x)\}](t)\Big\}\Big](u)=\bigg[\mathcal{I}\bigg\{\int_{-\infty}^{\infty}f(x)e^{-2\pi ixt}\ dx\bigg\}\bigg](u)$$
$$=\int_{-\infty}^{\infty}f(x)e^{-2\pi ix(u/i)}\ dx=\int_{-\infty}^{\infty}f(x)e^{-2\pi xu}\ dx=[\mathcal{K}\{f(x)\}](u)$$
Since $\mathcal{I}\mathcal{J}=\mathcal{K}$, $\mathcal{I}\mathcal{J}\mathcal{K}=\mathcal{K}\mathcal{K}=\mathcal{K}^2=\mathcal{P}$. As a result, the equation \\ $i^2=j^2=k^2=ijk=-1$ has an analogue for these operators:
$$\mathcal{I}^2=\mathcal{J}^2=\mathcal{K}^2=\mathcal{IJK}=\mathcal{P}$$
Using similar methods, we can construct the entire composition table, where row $\mathcal{I}$ and column $\mathcal{J}$ represents the operator $\mathcal{I}\mathcal{J}$. It is important to note that while the operators are applied from right to left, their corresponding quaternions are multiplied from left to right.
\begin{table}[H]
\centering
\begin{tabular}{|c?c|c|c|c|}
\hline
$\circ$ & $\mathit{1}$ & $\mathcal{I}$ & $\mathcal{J}$ & $\mathcal{K}$ \\ \Xhline{4\arrayrulewidth}
$\mathit{1}$ & $\mathit{1}$ & $\mathcal{I}$ & $\mathcal{J}$ & $\mathcal{K}$ \\ \hline
$\mathcal{I}$ & $\mathcal{I}$ & $\mathcal{P}$ & $\mathcal{K}$ & $\mathcal{P}\mathcal{J}$ \\ \hline
$\mathcal{J}$ & $\mathcal{J}$ & $\mathcal{P}\mathcal{K}$ & $\mathcal{P}$ & $\mathcal{I}$ \\ \hline
$\mathcal{K}$ & $\mathcal{K}$ & $\mathcal{J}$ & $\mathcal{P}\mathcal{I}$ & $\mathcal{P}$ \\ \hline
\end{tabular}
\end{table}
\noindent Also, although the operators themselves correspond to quaternions, the domain of both the original and transformed functions is $\mathbb{C}$. The same is true for the rest of the operators we will derive in this section. \\\null\\
The next step is to generalize $\mathcal{I}$, $\mathcal{J}$, and $\mathcal{K}$ to fractional values. This is also done using our matrix representation $a+bi+cj+dk\rightarrow\big[\begin{smallmatrix}
a+bi & c+di \\
-c+di & a-bi
\end{smallmatrix}\big]$. However, we now let $a=\cos\theta$ and $b$, $c$, or $d=\sin\theta$:
$$\mathrm{I}_{\theta}=\begin{bmatrix}
e^{i\theta} & 0 \\
0 & e^{-i\theta}
\end{bmatrix},\quad
\mathrm{J}_{\theta}=\begin{bmatrix}
\cos\theta & \sin\theta \\
-\sin\theta & \cos\theta
\end{bmatrix},\quad
\mathrm{K}_{\theta}=\begin{bmatrix}
\cos\theta & i\sin\theta \\
i\sin\theta & \cos\theta
\end{bmatrix}$$
Since these matrices still have determinant 1, each one still has an associated LCT. For a given operator $\mathcal{T}_{\theta}$, we will choose the constant factor appropriately to ensure that $\mathcal{T}_{n\pi/2}=\mathcal{T}^n\ \forall n\in\mathbb{Z}$.
$$[\mathcal{I}_{\theta}\{f(x)\}](u)=f(ue^{-i\theta})$$
$$[\mathcal{J}_{\theta}\{f(x)\}](u)=\sqrt{1-i\cot(\theta)}\int_{-\infty}^{\infty}e^{\pi i\cot(\theta)(x^2+u^2)}f(x)e^{-2\pi i\csc(\theta)xu}\ dx$$
$$[\mathcal{K}_{\theta}\{f(x)\}](u)=\sqrt{1-\cot(\theta)}\int_{-\infty}^{\infty}e^{\pi\cot(\theta)(x^2+u^2)}f(x)e^{-2\pi\csc(\theta)xu}\ dx$$
These operators correspond to the quaternions $e^{i\theta}$, $e^{j\theta}$, and $e^{k\theta}$. $\mathcal{I}_{\theta}$ is a scaling of the domain by a factor of $e^{i\theta}$, and $\mathcal{J}_{\theta}$ is the Fractional Fourier Transform. Likewise, $\mathcal{K}_{\theta}$ may be interpreted as the Fractional Laplace Transform. \\\null\\
The final step is to generalize $\mathcal{I}_{\theta}$, $\mathcal{J}_{\theta}$, and $\mathcal{K}_{\theta}$ to a single operator corresponding to any unit quaternion, or versor. To do this, we use a version of Hopf coordinates:
$$a=\cos\xi_{1}\cos\eta,\quad b=\sin\xi_{1}\cos\eta,\quad c=\cos\xi_{2}\sin\eta,\quad d=\sin\xi_{2}\sin\eta$$
$(\xi_{1},\eta,\xi_{2})$ are a set of angles that determine any point on the 3-sphere \\ $a^2+b^2+c^2+d^2=1$. This corresponds to the versor
$$a+bi+cj+dk=\cos\xi_{1}\cos\eta+i\sin\xi_{1}\cos\eta+j\cos\xi_{2}\sin\eta+k\sin\xi_{2}\sin\eta$$
Using the matrix representation one last time, we obtain
$$\mathrm{V}_{(\xi_{1},\eta,\xi_{2})}=\begin{bmatrix}
e^{i\xi_{1}}\cos\eta & e^{i\xi_{2}}\sin\eta \\
-e^{-i\xi_{2}}\sin\eta & e^{-i\xi_{1}}\cos\eta
\end{bmatrix}$$
Finally, we can find its associated LCT, choosing the constant factors appropriately so that $\mathcal{I}_{\theta}$, $\mathcal{J}_{\theta}$, and $\mathcal{K}_{\theta}$ are special cases. We will call this the Versor Transform:
\begin{adjustwidth}{-73pt}{-73pt}
$$[\mathcal{V}_{(\xi_{1},\eta,\xi_{2})}\{f(x)\}](u)=
\begin{cases}
\displaystyle\sqrt{1-ie^{-i\xi_{2}}\cot(\eta)}\int_{-\infty}^{\infty}e^{\pi ie^{-i\xi_{2}}\cot(\eta)(e^{i\xi_{1}}x^2+e^{-i\xi_{1}}u^2)}f(x)e^{-2\pi ie^{-i\xi_{2}}\csc(\eta)xu}\ dx, & \sin\eta\neq0 \\ \\
f(ue^{-i\xi_{1}}\cos(\eta)), & \sin\eta=0
\end{cases}$$
\end{adjustwidth}
$$\text{Notice that } \mathcal{V}_{(\theta,0,\xi_{2})}=\mathcal{I}_{\theta},\quad \mathcal{V}_{(0,\theta,0)}=\mathcal{J}_{\theta},\quad \mathcal{V}_{(0,\theta,\pi/2)}=\mathcal{K}_{\theta}.$$
The inverse transform $\mathcal{V}^{-1}_{(\xi_{1},\eta,\xi_{2})}$ is equivalent to $\mathcal{V}_{(-\xi_{1},-\eta,\xi_{2})}$. This can be seen by the relationship between their underlying matrices:
\begin{adjustwidth}{-25pt}{-25pt}
$$\mathrm{V}^{-1}_{(\xi_{1},\eta,\xi_{2})}=\begin{bmatrix}
e^{-i\xi_{1}}\cos\eta & -e^{i\xi_{2}}\sin\eta \\
e^{-i\xi_{2}}\sin\eta & e^{i\xi_{1}}\cos\eta
\end{bmatrix}=
\begin{bmatrix}
e^{i(-\xi_{1})}\cos(-\eta) & e^{i\xi_{2}}\sin(-\eta) \\
-e^{-i\xi_{2}}\sin(-\eta) & e^{-i(-\xi_{1})}\cos(-\eta)
\end{bmatrix}=\mathrm{V}_{(-\xi_{1},-\eta,\xi_{2})}$$
\end{adjustwidth}
Just as every unit complex number $e^{i\theta}$ had an associated Fractional Fourier Transform $\mathcal{F}_{\theta}$, every unit quaternion $e^{i\xi_{1}}\cos\eta+e^{i\xi_{2}}j\sin\eta$ has an associated Versor Transform $\mathcal{V}_{(\xi_{1},\eta,\xi_{2})}$. For example, any unit quaternion with zero real part is a square root of -1. Analogously, any Versor Transform $\mathcal{V}=\mathcal{V}_{(\xi_{1},\eta,\xi_{2})}$ with $\cos\xi_{1}\cos\eta=0$ will satisfy $\mathcal{V}^2=\mathcal{P}$, $\mathcal{V}^3=\mathcal{V}^{-1}$, and $\mathcal{V}^4=\mathit{1}$. \\\null\\
Let us examine one particular case of the Versor Transform. Define $\mathcal{H}_{\theta}$ as follows:
$$[\mathcal{H}_{\theta}\{f(x)\}](u)=[\mathcal{V}_{(\xi_{1},\pi/2,\theta)}\{f(x)\}](u)=\int_{-\infty}^{\infty}f(x)e^{-2\pi ie^{-i\theta}xu}\ dx$$
Then, $\mathcal{H}_{\theta}$ can be thought of as a hybrid between the Fourier and Laplace Transforms. For example, here are its values when $\theta$ is a multiple of $\frac{\pi}{2}$:
\begin{table}[H]
\centering
\begin{tabular}{|c?c|c|c|c|c|}
\hline
$\theta$ & $0$ & $\pi/2$ & $\pi$ & $3\pi/2$ & $2\pi$ \\ \hline
$\mathcal{H}_{\theta}$ & $\mathcal{J}$ & $\mathcal{K}$ & $\mathcal{J}^{-1}$ & $\mathcal{K}^{-1}$ & $\mathcal{J}$ \\ \hline
\end{tabular}
\end{table}
\noindent This encourages an alternate view of the relationship between the two transforms. In summary: \\\null\\
$\mathcal{I}$ takes a function from the real time domain into the imaginary time domain. $\mathcal{J}$ takes a function from the real time domain into the real frequency domain. $\mathcal{K}$ takes a function from the real time domain into the imaginary frequency domain. As a result, the domain of the Versor Transform can be considered a linear combination of both real and imaginary time and frequency.
\subsection{Identities}
These can be derived by letting
$$\mathrm{M}=\mathrm{V}_{(\xi_{1},\eta,\xi_{2})}=\begin{bmatrix}
e^{i\xi_{1}}\cos\eta & e^{i\xi_{2}}\sin\eta \\
-e^{-i\xi_{2}}\sin\eta & e^{-i\xi_{1}}\cos\eta
\end{bmatrix}$$
in the corresponding identities for the LCT. For the Convolution and Cross-Correlation Theorems, the constant factors are adjusted appropriately. Assume $\sin\eta\neq0$.
\subsubsection{Shift Identities}
For a shift in the original function's domain,
\begin{adjustwidth}{-49pt}{-49pt}
$$[\mathcal{V}_{(\xi_{1},\eta,\xi_{2})}\{f(x-x_{0})\}](u)=e^{\pi ie^{i\xi_{1}}e^{-i\xi_{2}}\cos(\eta)\sin(\eta)x_{0}^2}e^{-2\pi ie^{-i\xi_{2}}\sin(\eta)x_{0}u}[\mathcal{V}_{(\xi_{1},\eta,\xi_{2})}\{f(x)\}](u-e^{i\xi_{1}}\cos(\eta)x_{0})$$
\end{adjustwidth}
For a shift in the transformed function's domain,
\begin{adjustwidth}{-56pt}{-56pt}
$$e^{\pi ie^{-i\xi_{1}}e^{-i\xi_{2}}\cot(\eta)u_{0}^2}e^{-2\pi ie^{-i\xi_{1}}e^{-i\xi_{2}}\cot(\eta)uu_{0}}[\mathcal{V}_{(\xi_{1},\eta,\xi_{2})}\{e^{2\pi ie^{-i\xi_{2}}\csc(\eta)xu_{0}}f(x)\}](u)=[\mathcal{V}_{(\xi_{1},\eta,\xi_{2})}\{f(x)\}](u-u_{0})$$
\end{adjustwidth}
\subsubsection{Derivative Identities}
For the $n^{\text{th}}$ derivative of the original function,
\begin{adjustwidth}{-116pt}{-116pt}
$$\bigg[\mathcal{V}_{(\xi_{1},\eta,\xi_{2})}\bigg\{\frac{d^n}{dx^n}f(x)\bigg\}\bigg](u)=(-1)^n\sum_{k=0}^{n}(-2\pi ie^{-i\xi_{2}}\csc(\eta)u)^k\binom{n}{k}\bigg[\mathcal{V}_{(\xi_{1},\eta,\xi_{2})}\bigg\{e^{-\pi ie^{i\xi_{1}}e^{-i\xi_{2}}\cot(\eta)x^2}\frac{d^{n-k}}{dx^{n-k}}(e^{\pi ie^{i\xi_{1}}e^{-i\xi_{2}}\cot(\eta)x^2})f(x)\bigg\}\bigg](u)$$
\end{adjustwidth}
For the $n^{\text{th}}$ derivative of the transformed function,
\begin{adjustwidth}{-104pt}{-104pt}
$$\sum_{k=0}^{n}\binom{n}{k}e^{-\pi ie^{-i\xi_{1}}e^{-i\xi_{2}}\cot(\eta)u^2}\frac{d^{n-k}}{du^{n-k}}(e^{\pi ie^{-i\xi_{1}}e^{-i\xi_{2}}\cot(\eta)u^2})[\mathcal{V}_{(\xi_{1},\eta,\xi_{2})}\{(-2\pi ie^{-i\xi_{2}}\csc(\eta)x)^kf(x)\}](u)=\frac{d^n}{du^n}\Big([\mathcal{V}_{(\xi_{1},\eta,\xi_{2})}\{f(x)\}](u)\Big)$$
\end{adjustwidth}
\subsubsection{Parseval's Theorem}
For any two functions $f(x)$ and $g(x)$ for which their Versor Transforms exist:
$$\int_{-\infty}^{\infty}f(x)\overline{g(x)}\ dx=\int_{-\infty}^{\infty}[\mathcal{V}_{(\xi_{1},\eta,\xi_{2})}\{f(x)\}](u)\ \overline{[\mathcal{V}_{(\xi_{1},\eta,\xi_{2})}\{g(x)\}](u)}\ du$$
Letting $f(x)=g(x)$ gives Plancherel's Theorem:
$$\int_{-\infty}^{\infty}|f(x)|^2\ dx=\int_{-\infty}^{\infty}\Big|[\mathcal{V}_{(\xi_{1},\eta,\xi_{2})}\{f(x)\}](u)\Big|^2\ du$$
\subsubsection{Convolution Theorem}
\begin{adjustwidth}{-93pt}{-93pt}
$$\sqrt{1-ie^{-i\xi_{2}}\cot(\eta)}\ e^{\pi ie^{-i\xi_{1}}e^{-i\xi_{2}}\cot(\eta)u^2}[\mathcal{V}_{(\xi_{1},\eta,\xi_{2})}\{e^{-\pi ie^{i\xi_{1}}e^{-i\xi_{2}}\cot(\eta)x^2}[fh*gh](x)\}](u)=[\mathcal{V}_{(\xi_{1},\eta,\xi_{2})}\{f(x)\}\cdot\mathcal{V}_{(\xi_{1},\eta,\xi_{2})}\{g(x)\}](u)$$
\end{adjustwidth}
\begin{adjustwidth}{-97pt}{-97pt}
$$\sqrt{1+ie^{-i\xi_{2}}\cot(\eta)}^{-1}e^{-\pi ie^{-i\xi_{1}}e^{-i\xi_{2}}\cot(\eta)u^2}[\mathcal{V}_{(\xi_{1},\eta,\xi_{2})}\{e^{\pi ie^{i\xi_{1}}e^{-i\xi_{2}}\cot(\eta)x^2}[f\cdot g](x)\}](u)=[\mathcal{V}_{(\xi_{1},\eta,\xi_{2})}\{f(x)\}k*\mathcal{V}_{(\xi_{1},\eta,\xi_{2})}\{g(x)\}k](u)$$
\end{adjustwidth}
where $h(x)=e^{\pi ie^{i\xi_{1}}e^{-i\xi_{2}}\cot(\eta)x^2}$ and $k(u)=e^{-\pi ie^{-i\xi_{1}}e^{-i\xi_{2}}\cot(\eta)u^2}$.
\subsubsection{Cross-Correlation Theorem}
\begin{adjustwidth}{-96pt}{-96pt}
$$\overline{\sqrt{1-ie^{-i\xi_{2}}\cot(\eta)}}\ e^{-\pi ie^{-i\xi_{1}}e^{-i\xi_{2}}\cot(\eta)u^2}[\mathcal{V}_{(\xi_{1},\eta,\xi_{2})}\{e^{-\pi ie^{i\xi_{1}}e^{-i\xi_{2}}\cot(\eta)x^2}[fh\star gh](x)\}](u)=[\overline{\mathcal{V}_{(\xi_{1},\eta,\xi_{2})}\{f(x)\}}\cdot\mathcal{V}_{(\xi_{1},\eta,\xi_{2})}\{g(x)\}](u)$$
\end{adjustwidth}
\begin{adjustwidth}{-100pt}{-100pt}
$$\overline{\sqrt{1+ie^{-i\xi_{2}}\cot(\eta)}}^{-1}e^{-\pi ie^{-i\xi_{1}}e^{-i\xi_{2}}\cot(\eta)u^2}[\mathcal{V}_{(\xi_{1},\eta,\xi_{2})}\{e^{-\pi ie^{i\xi_{1}}e^{-i\xi_{2}}\cot(\eta)x^2}[\overline{f}\cdot g](x)\}](u)=[\mathcal{V}_{(\xi_{1},\eta,\xi_{2})}\{f(x)\}k\star\mathcal{V}_{(\xi_{1},\eta,\xi_{2})}\{g(x)\}k](u)$$
\end{adjustwidth}
where $h(x)=e^{\pi ie^{i\xi_{1}}e^{-i\xi_{2}}\cot(\eta)x^2}$ and $k(u)=e^{-\pi ie^{-i\xi_{1}}e^{-i\xi_{2}}\cot(\eta)u^2}$.
\newpage
\bibliographystyle{plain}
\bibliography{references}
\end{document}